%% file: neutro.tex
\documentclass[11pt,a4paper,final]{article}

\usepackage{url}
\usepackage[a4paper,portrait]{geometry}
\usepackage{amsmath,amsfonts,latexsym,amssymb,amsthm}
\usepackage{graphicx,color,ifpdf}

\input{macros}

\theoremstyle{plain}
\newtheorem{thm}{Theorem}[section]

\newtheorem{lem}[thm]{Lemma}
\theoremstyle{definition}

\theoremstyle{plain} % {remark}
\newtheorem{rem}[thm]{Remark}
\newtheorem{xpl}[thm]{Example}

\title{Explicit formulas for a continuous stochastic maturation model.\\
Application to anticancer drug pharmacokinetics/pharmacodynamics.}

\author{Djalil~\textsc{Chafa\"\i} and Didier~\textsc{Concordet}}

\date{March 2008.}

\begin{document}

%- Local macros.
%-.

\maketitle

\begin{abstract}
  We present a continuous time model of maturation and survival, obtained as
  the limit of a compartmental evolution model when the number of compartments
  tends to infinity. We establish in particular an explicit formula for the
  law of the system output under inhomogeneous killing and when the input
  follows a time-inhomogeneous Poisson process. This approach allows the
  discussion of identifiability issues which are of difficult access for
  finite compartmental models. The article ends up with an example of
  application for anticancer drug pharmacokinetics/pharmacodynamics.
\end{abstract}

{
\footnotesize
 \noindent%
 \textbf{Note: this article is accepted for publication in Stochastic Models 
  \copyright Taylor \& Francis, April 2008.} \\
 \noindent%
 \textbf{AMS-MSC-2000}: %
 92B05 % General biology and biomathematics %
 ; 60G07 % General theory of stochastic processes%
 ; 60K25 % Queueing theory %
 ; 60K20 % Applications of Markov renewal processes %
 ; 60J80 % Branching processes (Galton-Watson, birth-and-death, etc.)
 ; 60J85 % Applications of branching processes %
 ; 47D08 % Schrödinger and Feynman-Kac semigroups %
 ; 60J25 % Markov processes with continuous parameter %
 ; 60J75 % Jump processes %
 .%
 %\tableofcontents
}

%%%
\addcontentsline{toc}{section}{Introduction}%
\section*{Introduction}
%%%

Compartmental models are widely used in biology and medicine for the modelling
of evolution phenomena. In particular, these models are very usual in the
pharmacokinetics/pharmacodynamics analysis of toxic effects of anticancer drugs
on white blood cells. Several mathematical and computational aspects of such
models remain untouched. For instance, at least for catenary chain models,
there is no rigorous rule for the choice of the number of compartments, and
many identifiability issues are not elucidated. In the present study, we
consider the limit of a discrete family of compartmental models when the
number of compartments tends to infinity. The obtained limit looks simpler. It
provides explicit formulas for the mean occupation of the compartment of
interest, and allows the discussion of important identifiability issues.
Roughly speaking, we replace a finite catenary chain of compartments with
time-inhomogeneous rates by a two compartments time-inhomogeneous model with
lag. Much of the article is devoted to the mathematical derivation and
analysis of these models. We hope that beyond the mathematical aspects, the
main results -- illustrated on a simple example at the end of the article --
may serve the quantitative biologists.

%This study is in the spirit of previous works of Jacquez \& Simon and
%Schuhmacher \& Thieme. 
The mean occupation of the compartment of interest appears as the expectation
of an underlying stochastic process. This process is a time-inhomogeneous
M/M/$\infty$ queue, with an explicit biological interpretation of its input
and output rates. This leads to a nice explicit Binomial-Poisson formula for
the instantaneous law. In a way, our approach can be seen as a complement and
extension of the boxcartrain models considered for example in
\cite{MR860975,MR953566}, see also \cite{MR1013208,MR1950761} and references
therein. It is also closely linked with ``binomial catastrophe models'', see
for instance \cite{MR2052588} and references therein. The novelty is mainly
the space-time inhomogeneous killing, the stochastic interpretation in terms
of Feynman-Kac's formul\ae, and the computation of the output occupation law
when the input follows a time-inhomogeneous Poisson process. In practice, the
formulas that we derive provide a good compromise between computer time and
numerical errors, without loss in interpretation, as suggested by the example
presented at the end of the article. However,
pharmacokinetics/pharmacodynamics populational aspects are outside the scope
of this work. %, and will make hopefully the matter of a forthcoming article.

\textbf{Outline of the rest of the article.} 
The rest of the article is organized as follows.

In Section \ref{se:finite-compart}, we briefly introduce the
time-inhomogeneous deterministic compartmental systems with linear rates. In
particular, we present two stochastic interpretations in terms of particles.
The first one is based on interacting time-inhomogeneous M/M/$\infty$ queueing
systems, whereas the second is based on the occupation of independent random
walks on the graph of compartments, subject to birth and death. These
interpretations are at the heart of the results of Section
\ref{se:cont-sto-mat-mod} regarding limits of catenary chains of compartments.

In section \ref{se:cont-sto-mat-mod}, we present a time-inhomogeneous
maturation and survival model, and the related counting processes. We show how
a finite catenary chain of compartments leads, when the number of compartments
tends to infinity, to a simple time-inhomogeneous two compartments model with
lag. The main result is given by Theorem \ref{th:law}, which provides simple
explicit integral formulas for the compartment of interest. These formulas can
be seen in turn as a byproduct of Feynman-Kac's formulas, and allow the
discussion of identifiability issues in Section \ref{ss:identif}.

In Section \ref{se:example}, the results of Sections \ref{se:finite-compart}
and \ref{se:cont-sto-mat-mod} are illustrated on a simple example related to
anticancer drug toxicity pharmacokinetics/pharmacodynamics context. Namely, we
provide a comparison between classical finite catenary chains of compartments
in one hand, and our two compartmental limit with lag in the other hand.

%\textbf{Notations and conventions.} In the sequel, we denote by $\dR_+$ the
%non negative half line $[0,+\infty)$. For a set $A$, we denote by $\rI_A$ the
%indicator function defined by $\rI_A(x)=1$ if $x\in A$ and $\rI_A(x)=0$
%otherwise. We denote by $\cL(X)$ the law of the random variable $X$, by
%$\cL(X\,\vert\,Y)$ the conditional law of $X$ given $Y$, and by
%$\bE(X\,\vert\,Y)$ the conditional expectation of $X$ given $Y$. We use in
%addition the classical associated maps $\cL(X\vert Y=y)$ and $\bE(X\vert
%Y=y)$. We denote by $\cP(\la)$, the Poisson distribution with mean
%$\la\in\dR_+$, by $\cE(\la)$ the exponential distribution of mean $1/\la$, and
%by $\cB(n,p)$ the binomial distribution of size $n\in\dN$ and parameter
%$p\in[0,1]$. In particular, $\cB(1,p)$ is the Bernoulli distribution
%$p\de_1+(1-p)\de_0$ and $\cB(n,p)=\cB(1,p)^{\ast n}$ where $\ast$ denotes the
%convolution.

%%
\section{Finite compartmental systems with time-dependent linear rates}
\label{se:finite-compart}

Consider a system of compartments indexed by the finite set $I$. Each
compartment contains a quantity of matter. As we will see in the sequel, the
amount of matter can be represented by a discrete number, or by a real number,
depending on the interpretation chosen. The matter is the subject of transfers
between compartments. It can also be created (external inflow) or destroyed
(outflow) in each compartment.

The reader may find an accessible introduction to simple compartmental systems
with several examples in the recent book \cite{MR1765331} by Matis and Kiffe,
see also \cite{MR1207797} and references therein. The study of the total
amount of matter in the system is for instance addressed in \cite{MR1950752},
see also \cite{MR1202398}. The reader may find a study of time-delayed
compartmental systems in \cite{MR2092594} and \cite{MR1950761} and references
therein. The literature regarding the stochastic interpretations of
time-inhomogeneous compartmental systems is less rich. Let us recall the
essential aspects of finite compartmental systems with linear rates.

\subsection{Deterministic systems and linear ordinary differential equations}

\begin{figure}[hbpt]
  \begin{center}
    \ifpdf\input{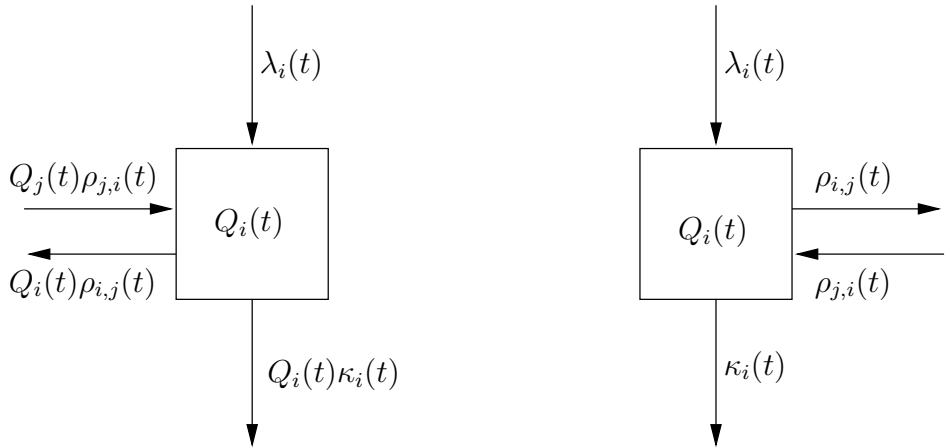}\else\input{compgen.pstex_t}\fi
    \caption{The left hand side diagram shows the flows in a compartmental
      system. Here, the transfer rates and the outflow depends linearly on
      the content of the compartment that gives matter. In contrast, the
      inflow does not depend on the content of the compartments. It is
      customary to make the $Q_i(t)$ and $Q_j(t)$ parts of the rates implicit,
      as show in the right hand side diagram.}
    \label{fi:compgen}
  \end{center}
\end{figure}

We consider here that the amount of matter is represented by a non-negative
real number. Let $Q_i(t)$ be the amount of matter in compartment $i$ at time
$t\geq t_0$, where $t_0$ is the initial time. We denote by $Q(t)$ the vector
$i\in I\mapsto Q_i(t)$, and we identify the set $\dR^I$ with $\dR^n$ where
$n:=\mathrm{card}(I)$. Let us consider the dynamics of $(Q(t))_{t\geq t_0}$
described by the system of linear ordinary differential equations
\begin{equation}\label{eq:determ-syst}
  \forall i\in I,\ \forall t\geq t_0,\ \pd_t Q_i(t) %
  = \la_i(t) %
  + \sum_{j\neq i} \rho_{j,i}(t) Q_j(t) %
  - Q_i(t)\SBRA{\kappa_i(t) %
    + \sum_{j\neq i} \rho_{i,j}(t)},
\end{equation}
with initial condition $Q(t_0)$. Such a dynamics is described by figure
\ref{fi:compgen}. Here $\la_i(t)$ is a creation rate for compartment $i$
(inflow), $\kappa_i(t)$ is the destruction rate for compartment $j$ (outflow),
and $\rho_{i,j}(t)$ with $i\neq j$ is the transfer rate from compartment $i$
to compartment $j$. These rates are non-negative. Apart from the $\la_i$
rates, the rates $\kappa_i$ and $\rho_{i,j}$ act proportionally to the content
of the compartment that gives matter. Note also that a $\kappa_i$ with
negative values might mimic a proportional auto-inflow. In vector/matrix
language, the dynamics of $Q(t)$ is of the form
\begin{equation}\label{eq:ode}
\pd_t Q(t) = \bM(t) Q(t) + \la(t),
\end{equation}
where $\la(t)$ denote the ``vector'' $(\la_i(t);i\in I)$, and where $\bM(t)$
is the matrix defined by
$$
\bM_{i,j}(t)
:=
\begin{cases}
 \rho_{j,i}(t) & \text{if $i\neq j$} \\
 -\sum_{k\neq i}\rho_{i,k}(t)-\kappa_i(t) & \text{if $i=j$}
\end{cases}
$$
for any $i,j\in I$. The traditional theory of linear ordinary differential
equations gives
$$
Q(t)=\bR(t,t_0)Q(t_0)+\int_{t_0}^t\!\bR(t,u)\la(u)\,du
$$
for any $t\geq t_0$, where the resolvent $\bR$ is the solution of the matrix
linear differential equation $\bR(t_0,t_0)=\bI$ and $\pd_t
\bR(t,t_0)=\bM(t)\bR(t,t_0)$ for any $t\geq t_0$. In the literature related to
compartmental systems, $\bM$ is sometimes referred as the ``transfer matrix''.
When $\bM$ \emph{does not depend on time}, $\bR(v,u)=\exp((v-u)\bM)$. This
matrix exponential can be explicitly computed in very special situations where
for example $\mathrm{card}(I)$ is small or where $\bM$ has a nice structure.
Usually, the method is to diagonalize $\bM$, which leads to an expression of
the solution as a linear combination of exponential functions related to the
spectrum of $\bM$. There is no simple closed formula for $Q(t)$ when the rates
depend on time.

\begin{figure}[hbpt]
  \begin{center}
  \ifpdf\input{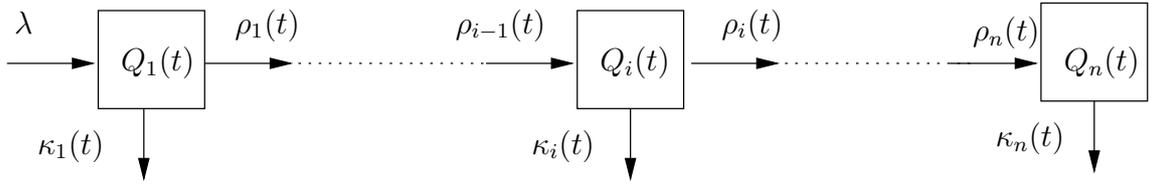}\else\input{fincat.pstex_t}\fi
  \caption{Finite catenary chain of $n$ compartments with single inflow at the
    chain head.}
  \label{fi:fincat}
  \end{center}
\end{figure}

\begin{xpl}[Finite catenary chain of compartments]\label{xp:finite-catenary}
  It corresponds to a finite system of $n$ compartments, labelled by
  $I=\{1,\ldots,n\}$, for which $\la_i\equiv0$ if $i\neq 1$, and
  $\rho_{i,j}\equiv0$ if $j\neq i+1$. It is customary to abridge $\la_1$ by
  $\la$ and $\rho_{i,i+1}$ by $\rho_i$. In such a compartmental system, there
  is only one external inflow with rate $\la$ at the left extremity of the
  chain. Moreover, the interaction of a compartment is limited to its right
  neighbor. The topology of this compartmental system is depicted in figure
  \ref{fi:fincat}. The matrix $\bM$ is band diagonal, which makes possible yet
  tedious the explicit computation of its exponential, when the rates do not
  depend on time.
\end{xpl}

\subsection{Stochastic interpretation}
\label{ss:stoch-interp}

Let us consider a finite system of $n$ compartments, labelled by $I$, and with
rates $\la$, $\kappa$, and $\rho$. In the stochastic interpretation, each
compartment contains particles, and we denote by $N_t^i$ the number of
particles in compartment $i\in I$ at time $t\geq t_0$. The particles are
indistinguishable, and can be created, destroyed, or can move from one
compartment to another, according to a Markovian dynamics of the vector
$N_t:=(N_t^i;i\in I)$. For a good choice of Markovian dynamics, the average
number $\dE(N_t\,\vert\,N_{t_0}=x)$ of particles per compartment is the
solution of a system of linear differential equations with initial condition
$x$. The reader may find a general presentation of Markov processes and
related topics in the book \cite{MR838085} by Ethier \& Kurtz. Kelly gave in
\cite[sec. 4.5 p. 113--117]{MR554920} an accessible introduction to the
stochastic interpretation of compartmental systems. Here we focus on the
time-inhomogeneous case.

\begin{thm}\label{th:ode-inhom}
  Let $I$ be a finite set. Let $(N_t)_{t\geq t_0}$ be a time-inhomogeneous
  Markov process with state space $E:=\dN^I$, and generators $(\bL_t)_{t\geq
    t_0}$. Assume in addition that for any $t\geq t_0$, the function $x\in
  E\mapsto A_t(x):=\sum_{y\in E}\bL_{t,x,y} y$ is affine. For any $x\in E$,
  and any $t\geq s\geq t_0$, let $Q(s,t,x):=\dE(N_t\,\vert\,N_s=x)$. Then $Q$
  is the solution of the linear differential equation
\begin{equation*}
\pd_t Q(s,t,x) = A_t(Q(s,t,x)),
\end{equation*}
for any $t\geq s\geq t_0$, with initial condition $Q(s,s,x)=x$.
\end{thm}

\begin{proof}
  The result is a consequence of Chapman-Kolmogorov equations. Namely, for any
  $t\geq s\geq t_0$ and any function $f:E\to\dR$, let $P_{s,t}(f)$ be the
  function $E\to\dR$ defined by $P_{s,t}(f)(x):=\dE(f(N_t)\,\vert\,N_s=x)$. In
  particular, $P_{s,s}(f)=f$. The Markovianity of the process is captured by
  the Chapman-Kolmogorov equation which writes
  $P_{s,u}(P_{u,t}(f))=P_{s,t}(f)$ for any $t\geq u\geq s\geq t_0$. Recall
  that $\bL_t$ acts linearly on $f$ as a matrix. We denote by $\bL_t f$ the
  function $E\to\dR$ defined by $(\bL_t f)(x):=\sum_{y\in E}\bL_{t,x,y}f(y)$.
  The definition of $\bL$ gives
  $$
  (\bL_t f)(x):=
  \lim_{\veps\to0^+}\frac{P_{t,t+\veps}(f)(x)-P_{t,t}(f)(x)}{\veps}.
  $$
  This yields, by using the Chapman-Kolmogorov equation, to the ``forward
  equation''
  \begin{equation*}
    \pd_t P_{s,t}(f)(x)
    :=\lim_{\veps\to0^+}\frac{P_{s,t+\veps}(f)(x)-P_{s,t}(f)(x)}{\veps}
    =P_{s,t}(\bL_t(f))(x).
  \end{equation*}
  The desired results follows immediately by considering the ad-hoc $f$
  function. Namely, for any $i\in I$, consider the function defined by
  $f(x)=x_i$ for any $x\in E$. One has $(\bL_t f)(x) = A_t(x)_i$. Since $A_t$
  in affine and $P_{s,t}$ is Markov, they commute and we get $P_t(\bL_t
  f)(x)=A_t(P_{s,t}(f)(x))_i$.
\end{proof}

Note that the Chapman-Kolmogorov equation which appears in the proof of
theorem \ref{th:ode-inhom} gives also the ``backward'' equation $\pd_s
P_{s,t}(f)(x)=-\bL_s(P_{s,t}(f))(x)$, which yields the system of linear
differential equations
\begin{equation*}%\label{eq:ode-inhom-1}
\pd_s Q(s,t,x) = -\sum_{y\in E}\bL_{s,x,y} Q(s,t,y)
\end{equation*}
with initial condition $Q(s,s,x)=x$. Of course, when the process is
time-homogeneous, $\bL_t$ does not depend on time anymore and
$P_{s,t}=P_{t_0,t_0+t-s}$. Thus, in that case, $P_u(\bL(f))=\bL(P_u(f))$ for
any $u\geq t_0$, which makes the forward and backward equations identical up
to a sign.

Let us give now an $\bL_t$ matrix on $E$ such that $A_t$ is exactly the affine
function of the deterministic compartmental system \eqref{eq:ode} constructed
from the triplet of rates $\la(t)$, $\kappa(t)$ and $\rho(t)$. We identify the
countable space $E:=\dN^I$ with $\dN^n$ where $n:=\mathrm{card}(I)$. We denote
by $e_1,\ldots,e_n$ the canonical basis of $\dR^n$, embedded in $E$. For any
$x,y\in E$ with $x\neq y$, we set
\begin{equation}\label{eq:def-Q-homo}
\bL_{t,x,y}
:=
\begin{cases}
  \la_i(t) & \text{if $y=x+e_i$ for $i\in I$ (birth of a particle)} \\
  x_i\rho_{i,j}(t) & \text{if $y=x-e_i+e_j$ for $i\neq j$ in $I$ (transfer of a particle)} \\
  x_i\kappa_i(t) & \text{if $y=x-e_i$ for $i\in I$ (death of a particle)} \\
  0 & \text{otherwise}
\end{cases}.
\end{equation}
The diagonal terms are such that the sum over each row is zero. The
transition related to $\rho$ can equivalently expressed as $x_j\rho_{j,i}(t)$
if $y=x-e_j+e_i$ for $j\neq i$. Recall that $N_t$ is an $n$-components vector
which represents the amount of particles in each of the $n$ compartments at
time $t$. A simple computation shows that for any $x\in E$ and any $i\in I$,
\begin{equation}\label{eq:prop-Q}
A_t(x)_i:=\sum_{y\in E}\bL_{t,x,y} y_i=
\la_i(t)
+\sum_{j\in I,j\neq i}x_j\rho_{j,i}(t)
-x_i\SBRA{\kappa_i(t)+\sum_{j\in I,j\neq i}\rho_{i,j}(t)}.
\end{equation}
We recognize immediately the equation \eqref{eq:ode} of the deterministic
compartmental system.

\subsubsection{Interpretation as interacting M/M/$\infty$ queues}

The generator \eqref{eq:def-Q-homo} suggests an interpretation of the process
$(N_t)_{t\geq t_0}$ as an $n$-dimensional M/M/$\infty$ queueing system with
interactions. In such a picture, each particle is a client, and each
compartment is an M/M/$\infty$ queue. The quantity $N_t^i$ is thus the number
of customers in the $i^\text{th}$ queue at time $t$. For each compartment
$i\in I$, and in absence of interaction ($\rho\equiv0$), the arrival rate is
$\la_i$, and the service rate is $\kappa_i$. Since the queue is M/M/$\infty$,
each new client gets immediately its own dedicated server, which explains the
coefficient $x_i\kappa_i$ in $\bL$. The interaction $\rho$ allows the clients
to move from queue $i$ to queue $j$ with rate $x_i\rho_{i,j}$. Thus, in
presence of interactions, the arrival rate in queue $i$ is $\la_i+\sum_j
x_j\rho_{j,i}$ and the service rate is $x_i\kappa_i+x_i\sum_j \rho_{i,j}$. The
clients in the queues are not ordered since these queues are of M/M/$\infty$
type. The random vector $N_t$ gives the number of clients in each queue at
time $t$. The dynamics of $(N_t)_{t\geq t_0}$ can be thus interpreted as $n$
interacting M/M/$\infty$ queues, e.g. as $n$ interacting birth and death
processes on $\dN$. These interacting queues can be seen as special Markovian
interacting particle systems, see \cite{MR1383122,MR2108619}. In statistical
mechanics, they constitute a ``zero-range dynamics on the graph of
compartments'', see \cite{MR2322692}. For the reader familiar with queueing
systems, these interacting queues can be seen as time-inhomogeneous ``Jackson
networks'' \cite{MR0093061,jackson,MR1996883}. These queues are ``in tandem''
in the case of a catenary chain of compartments with killing only at the end
of the chain.

In particular, with such an interpretation in mind, a one compartment system
is equivalent to a single M/M/$\infty$ queue. More generally, a catenary chain
of $n$ compartments is equivalent to a system of $n$ interacting M/M/$\infty$
queues, in which only the first one has an external inflow, the others being
inflowed by the clients that leave the preceding queue due to the interaction.

\begin{figure}[hbpt]
  \begin{center}
    \ifpdf\input{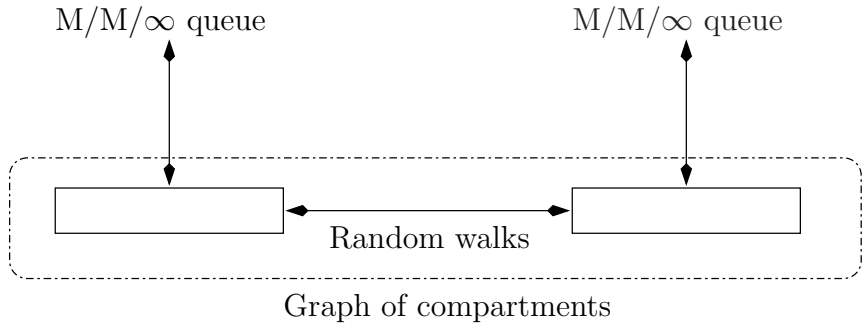}\else\input{stoch.pstex_t}\fi
    \caption{The two stochastic interpretations of finite compartmental systems
      with linear rates. Interacting M/M/$\infty$ queues versus independent
      random walks on the graph of compartments subject to birth and death.}
    \label{fi:stochint}
  \end{center}
\end{figure}

\subsubsection{Interpretation as independent particles on the graph of compartments}

Following \cite[sec. 4.5 p. 113--117]{MR554920}, by removing the
indistingishability of the particles assumed in the preceding interpretation,
$N_t$ is the occupation vector of the compartments at time $t$ for a
superposition of independent particles moving in $I$. These particles are
subject to birth ($\la$ rate) and death ($\kappa$ rate), in addition to their
motion in the system of compartments viewed as an oriented graph with $n$
vertices. This graph is a subset of $I\times I$, and $\rho_{i,j}\equiv0$ means
that the edge $(i,j)$ does not exist in the graph. In other words, $\rho$
defines an oriented graph, i.e. a topology on the system of compartments.

Note that when $n=1$ (trivial graph), we just have constructed the
M/M/$\infty$ queue from independent particles with survival rate $\kappa$,
arriving according to an independent Poisson point process of intensity $\la$.
When $n\geq 1$ and $\rho\equiv0$, the process consists in $n$ independent
M/M/$\infty$ queues.

In this interpretation, the creation of particles on the $n$ sites is obtained
with $n$ independent Poisson point processes, one for each external inflow.
The independence of the particles moving on the graph is a consequence of the
linearity in $x$ of the coefficients related to the $\rho$ rates in
$\bL_{x,y}$. The single particle dynamics involves the $\rho$ rates for its
motions. The $\kappa$ rates are incorporated by the addition of an external
cemetery state to the system, whereas the $\la$ rates are incorporated by
using an independent Poisson point process and independent copies of the
particle. It is known that adding an external cemetery state leads to
Feynman-Kac type formulas related to a very simple potential, see
\cite{MR2044973} for instance.

We have thus two stochastic interpretations of the dynamics
\eqref{eq:def-Q-homo} of $(N_t)$, as presented by figure \ref{fi:stochint}.
The first one is ``vertical'' and corresponds to interacting birth and death
processes in $\dN$ (M/M/$\infty$ queues subject to interactions), whereas the
second one is ``horizontal'' and corresponds to independent particles moving
on the graph of compartments (random walks on $I$ subject to birth and death).
Both hold for time-dependent rates.

\section{A continuous stochastic maturation model}
\label{se:cont-sto-mat-mod}

In the sequel, we consider a time-inhomogeneous stochastic maturation model
with killing, which provides a Poisson-Binomial formula for the occupation law
of the maturation system. The probability of survival in such a model can be
for instance the probability of traversal of a finite catenary chain of
compartments. We show that this probability has a nice explicit formula when
the length of the chain tends to infinity. Such a result is in accordance with
the stochastic interpretation of finite compartmental systems presented in
Section \ref{se:finite-compart}, in terms of time-inhomogeneous M/M/$\infty$
queues.

\subsection{Counting mature particles under killing and Poisson production}
\label{se:count-mat-simple}

Consider a maturation system in which the maturation takes a deterministic
duration $\tau>0$. A particle which begins its maturation at time $T\in\dR$
achieves its maturation at time $T+\tau$. We say that $T+\tau$ is the
\emph{maturation time} of the particle. Assume in addition that the particle
can die during the maturation process, with probability $1-p(T)$, where
$p(T)\in[0,1]$ is a deterministic real number which may depend on $T$. In such
a maturation with killing scheme, a particle which begins its maturation at
time $T\in\dR$ can either die with probability $1-p(T)$, or survive with
probability $p(T)$ and achieve its maturation at time $T+\tau$. Now, consider
independent particles which begin their maturation at times $T_0,T_1,\ldots$
These particles survive with probabilities $p(T_0),p(T_1),\ldots$, and achieve
their maturation at times $M_0=T_0+\tau,M_1=T_1+\tau,\ldots$ The following
Lemma is a consequence of the stability by thinning of Poisson point
processes, see for instance \cite{MR1207584} or \cite{MR1463943}.

\begin{lem}\label{le:thinning}
  Suppose that $T_0,T_1,\ldots$ are random and distributed according to a
  Poisson process with intensity $\la:\dR\to\dR_+$. Assume that this random
  process is independent from the maturation process of the particles. The
  maturation times $M_0,M_1,\ldots$ follow a Poisson point process on $\dR$
  with intensity $\la^*:\dR\to\dR_+$ given by
  $$
  \la^*(t):=\la(t-\tau)p(t-\tau).
  $$
\end{lem}

In other words, the number of mature particles produced during a bounded time
interval $I\subset\dR$ follows a Poisson distribution with mean
$$
\int_I\!\la^*(u)\,du=\int_I\!\la(u-\tau)\,p(u-\tau)\,du.
$$
Moreover, for any disjoint and bounded time intervals $I$ and $J$, the number
$N_I$ and $N_J$ of mature particles produced during these time intervals are
independent Poisson random variables. The reader familiar with point processes
may notice that the resulting point process is thus an inhomogeneous marked
Poisson point process, with Bernoulli marks.

\begin{figure}[hbpt]
  \begin{center}
    \ifpdf\input{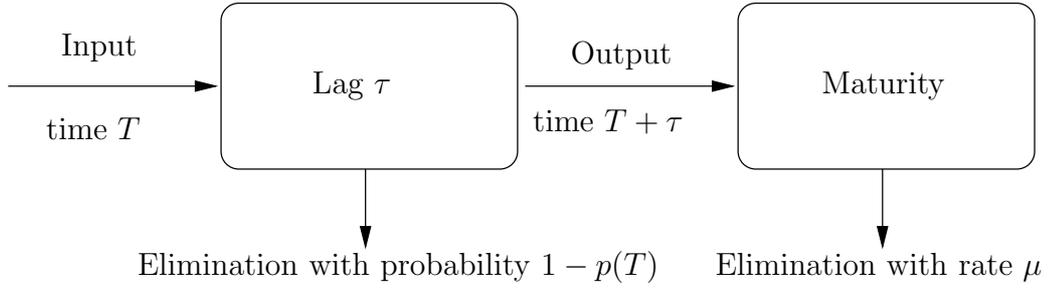}\else\input{schema2.pstex_t}\fi
    \caption{Schematic view of a maturation system with killing. Recall that
      $\tau$ is the maturation duration. The quantity $p(T)$ is the survival
      probability of a particle which has begun its maturation at time $T$.
      The quantity $\mu$ is the survival rate of a mature particle. The left
      hand side compartment (LC) corresponds to a Poisson deformation box with
      lag, whereas the right hand side compartment (RC) corresponds to an
      M/M/$\infty$-like counting process. The output of the LC is a point
      process which serves as an input for the RC.}
    \label{fi:schem}
  \end{center}
\end{figure}

Let us assume now in addition that independently, each mature particle
survives after maturation during a random duration with intensity
$\mu:\dR_+\to\dR_+$. In other words, given that a mature particle is still
alive at time $s$, and if $S$ is its remaining survival duration, then
$\dP(S>t)=\exp(-\int_{s}^{t}\!\mu(w)\,dw)$ for any $t\geq s$. The following
Theorem expresses the law of the number $N_t$ of mature particles still alive
at time $t\in\dR_+$. Figure \ref{fi:schem} gives a schematic view of the
system.

The following Theorem is a direct consequence of Lemma \ref{le:thinning}, but
can also be seen alternatively as a special case of a result of Keilson \&
Servi \cite{MR1274723} on time-inhomogeneous M/G/$\infty$ queuing processes
(see also \cite{foley,MR611802,taylor-boucherie} for related problems). Here
the resulting counting process $(N_t)_{t\geq0}$ is a time-inhomogeneous
M/M/$\infty$ queue. The reader may find a general introduction to Queuing
Processes in many books such as
\cite{MR1118840,MR1835969,MR1331145,MR554920,MR1207584,MR1996883}.

% , see for instance the reference books and \cite{MR1463943}.

\begin{thm}\label{th:maturgen}
  For any times $0\leq s\leq t$ and any $m\in\dN$,
  $$
  \cL(N_t\,\vert\,N_s=m) = \cB(m,\al(s,t))*\cP(\be(s,t))
  $$
  where
  $$
  \al(s,t):=\exp\PAR{-\int_s^t\!\mu(w)\,dw}
  \text{\quad and\quad}
  \be(s,t):=\int_s^t\!\la(u-\tau)p(u-\tau)\al(u,t)\,du.
  $$
  In particular, $\moy{}{N_t\,\vert\,N_s=m}=m\al(s,t)+\be(s,t)$ for any $0\leq
  s\leq t$ and any $m\in\dN$.
\end{thm}

Note that there are three sources of randomness here, which are supposed
independent: the initial time of maturation, the survival during the
maturation, and the survival duration after maturation.

\subsection{Modelling the maturation process itself for one particle}

Let us consider a single particle which begins its maturation at time
$T\in\dR$. This section is devoted to the modelling of the maturation and
survival process itself, by means of a Markovian dynamics on the time interval
$[T,+\infty)$. In particular, it provides an explicit formula of the
probability of survival during maturation $p(T)$.

\subsubsection{Finite state model}

\begin{figure}[hbpt] 
  \begin{center}
    \ifpdf\input{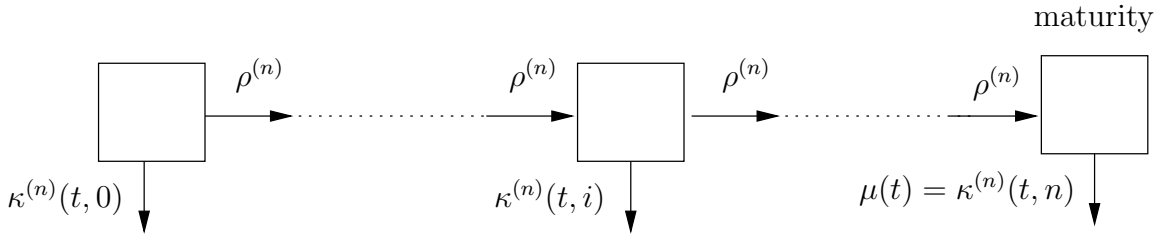}\else\input{schema0.pstex_t}\fi
    \caption{A catenary chain of compartments with one-way flow and possible
      killing, used as a maturation model maturation (or evolution) rate
      $\rho^{(n)}$, and killing rate $\kappa^{(n)}$. The down arrows lead to
      the cemetery state, which is not represented here.}
    \label{fi:family2}
  \end{center}
\end{figure}

We begin by considering a maturation model which consists in $n+1$ maturation
steps, with possible killing. More precisely, define the finite set
$$
\cS^{(n)}:=\{0,1,\ldots,n\}\cup\{c\}
$$
where $c\not\in\{0,1,\ldots,n\}$ will serve as a cemetery state. The
maturation of the particle is modelled by a motion from state $0$ to state
$n$. State $n$ corresponds to ``maturity''. The killing rate is modelled by a
function
\begin{equation*}
  \kappa^{(n)}:\dR\times\{0,\ldots,n\}\to\dR_+
\end{equation*}
which is smooth in the first variable (time). This killing rate may thus vary
in time and in space. The maturation rate is modelled by a constant positive
real number $\rho^{(n)}$. Consider now the continuous time Markov process
$(X^{(n)}_t)_{t>T}$ with state space $\cS^{(n)}$, and generator given for any
$i\neq j$ in $\cS^{(n)}$ and any $t\geq T$ by
\begin{equation}\label{eq:def-gi-discr-matrix}
(\bL^{(n)}_t)_{i,j}=
\begin{cases}
 \rho^{(n)} & \text{if $i\neq c$ and $j=i+1$} \\
 \kappa^{(n)}(t,i) & \text{if $i\neq c$ and $j=c$} \\
 0 & \text{otherwise}
\end{cases}.
\end{equation}
The diagonal terms of the matrix $\bL^{(n)}$ are such that the sum over each
row is $0$. The dynamics of the maturation process $(X^{(n)}_t)_{t\geq T}$
can be read directly on the expression of the generator above. Namely, at time
$t$ and from state $i$, the process can only move to state $i+1$ with rate
$\rho^{(n)}$, or die (being killed) with rate $\kappa(t,i)$ and placed in the
cemetery state $c$ from which it cannot escape. The maturation steps are
depicted by figure \ref{fi:family2}. In absence of killing, i.e. when
$\kappa^{(n)}\equiv0$, the process $(X^{(n)}_t)_{t\geq T}$ is a simple Poisson
process of intensity $\rho^{(n)}$, stopped when it reaches state $n$.

In such a maturation with killing model, we are interested in the mean
occupation of the state $n$ for a Poissonian number of such processes evolving
independently. One can define $p(T)$ as the probability of hitting state $n$
starting from state $0$:
$$
p(T)=\dP(\exists t>T; X^{(n)}_t=n\,\vert\,X^{(n)}_T=0)
\text{\quad and \quad} \mu(t)=\kappa^{(n)}(t,n).
$$
However, the corresponding maturation time is random, depends on $T$, and is
given by
$$
\tau=\inf\{t>T;X^{(n)}_t=n\}
$$
on the event $\{X^{(n)}_T=0,\exists t>T;X^{(n)}_t=n\}$. We will see in the
sequel that letting $n$ tends to $\infty$ provides a very simple explicit
formula for $p(T)$ and makes $\tau$ deterministic, without any loss in
interpretation.

\subsubsection{Continuous state limit model}

It is possible to express the law of $(X^{(n)})_{t\geq T}$ by mean of a
Feynman-Kac formula. Rather than using this approach, we prefer to consider
now the limit process obtained by letting the number of compartments $n$ tends
to $\infty$. We will use a scaling which corresponds to the convergence
$$
\BRA{\frac{0}{n},\frac{1}{n},\ldots,\frac{n}{n}} \to [0,1].
$$
Namely, assume that $\rho^{(n)}=n\rho$ for some positive constant $\rho$.
Assume that
$$
\lim_{n\to\infty}\sup_{i\in\{0,\ldots,n\};t\in\dR}%
\ABS{\kappa^{(n)}(t,i)-\kappa(t,i/n)}=0
$$
for some bounded piecewise continuous function
\begin{equation*}
  \kappa:\dR\times[0,1]\to\dR_+.
\end{equation*}

\begin{thm}\label{th:markov-lim}
  Let $x\in[0,1]$. On the event $\{X^{(n)}_T=\PENT{nx}\}$, and with the
  convention $c/n=c$, the rescaled process $(n^{-1}\,X_{t}^{(n)})_{t\geq T}$
  converges in law, when $n$ goes to $\infty$, toward the Markov process
  $(X_t)_{t\geq T}$ with state space $\cS:=[0,1]\cup\BRA{c}$ defined for any
  $T\leq s \leq t$ by
\begin{equation*}
  \cL\PAR{X_t\,\vert\,X_s=x} =
\begin{cases}
  p(s,t,x)\,\de_{x(t-s)}+(1-p(s,t,x))\,\de_c & \text{if $x\neq c$} \\
  \de_c & \text{if $x=c$}
\end{cases},
\end{equation*}
where $x(u):=\min(1,x+\rho(u-T))$, $p(u,v,x) := \exp\PAR{-\int_u^v
  \kappa(w,x(w-u))\,dw}$, and $p(u,v,c):=0$ for any $x\neq c$.
\end{thm}

\begin{proof}
  Consider the natural inclusion $\pi_n:\cS^{(n)}\to\cS$ defined by
  $\pi_n(i)=i/n$ for any $i\in\{0,\ldots,n\}$ and $\pi_n(c)=c$. By this way,
  any function $f:\cS\to\dR$ induces naturally a function $f(\pi_n)$ from
  $\cS^{(n)}$ to $\dR$, and one has
  \begin{equation*}
    \bL_t^{(n)}(f(\pi_n))(i)
    = n\rho\PAR{f\PAR{\frac{i+1}{n}}-f\PAR{\frac{i}{n}}}
    + \kappa\PAR{t,\frac{i}{n}}\PAR{f(c)-f\PAR{\frac{i}{n}}}.
  \end{equation*}
  By following step by step the classical argument based on Taylor's formulas
  presented in \cite[p. 1-5]{MR838085}, one can show that the Markov process
  $(n^{-1}\,X_t^{(n)})_{t\geq T}$ converges in law, when $n$ goes to $\infty$,
  toward the Markov process $(X_t)_{t\geq T}$ with state space
  $\cS:=[0,1]\cup\BRA{c}$ and infinitesimal generators
  \begin{equation} \label{eq:gi-matu}
    \bL_t(f)(x) =
    \begin{cases}
      \rho f'(x)+\kappa(t,x)\PAR{f(c)-f(x)} & \text{if $x\neq c$} \\
      0 & \text{if $x=c$}
    \end{cases}
  \end{equation}
  defined for continuous functions $f:\cS\to\dR$ which are smooth on $[0,1]$
  and vanish at the boundaries. %One can find detailed techniques of such
  %limiting procedures in \cite{MR838085} for example.
  The scaling method that we used is related to the \emph{fluid limits} of
  Queuing Theory (see for instance in \cite[ch. 9]{MR1996883}) and the
  \emph{hydrodynamic limits} of interacting particles systems in Statistical
  Mechanics (see for instance \cite{MR1707314}). The obtained generator
  $\bL_t$ above is the addition of a deterministic constant drift term $\rho
  f'(x)$ together with a random space-time inhomogeneous killing term
  $\kappa(t,x)\PAR{f(c)-f(x)}$. A Feynman-Kac formula or a direct check gives
  the desired expression of the law. Actually, the addition of a cemetery
  state corresponding to a killing term to a Markov process always gives rise
  to a Feynman-Kac formulae, see \cite{MR2044973,MR1121940}. In some ways, our
  case \eqref{eq:gi-matu} is the simplest one: killing term $+$ deterministic
  generator (constant drift).
  
\end{proof}

The interpretation in terms of the position $X_t$ of the particle is quite
clear: the particle is moving from left to right on the interval of positions
$[0,1]$ at constant speed $\rho$. The particle stops its motion when (and if)
it reaches the right extremity of $[0,1]$. At any time $t\in\dR_+$ and any
position $x\in[0,1]$, it can be killed (and thus placed in the cemetery state
$c$) with rate $\kappa(t,x)$. A particle is mature when it reaches the right
extremity $1$ of $[0,1]$. The positions $[0,1)$ correspond to the steps of
maturation. For any $t\in\dR$ and any $x\in[0,1]$, let us define
$g(t,x):=\kappa(t,x)\rI_{[0,1)}(x)$ and $\mu(t):=\kappa(t,1)$ in such a way
that
\begin{equation}\label{eq:def-k}
 \kappa(t,x)=g(t,x)\rI_{[0,1)}(x)+\mu(t)\rI_{\{1\}}(x).
\end{equation}
Function $g$ gathers the killing rate during maturation, whereas function
$\mu$ captures the killing rate after maturation. In the settings of Section
\ref{se:count-mat-simple}, the process $(X_t)_{t\geq T}$ corresponds to a
maturation model for which
\begin{equation}\label{eq:def-pt-mut}
\mu(t)=\kappa(t,1)
\text{\quad and\quad}
p(T)=p(T,T+\tau,0)=\exp\PAR{-\int_0^{\tau}\!g(T+w,\rho w)\,dw},
\end{equation}
where the ``maturation time'' $\tau$ is deterministic, independent of $t$, and
given by
\begin{equation}\label{eq:def-tau}
  \tau:=\frac{1}{\rho}.
\end{equation}
When the particle is not killed, i.e. on the event $\{X_T=0,\exists
s>0;X_s=1\}$, $\tau$ is exactly the deterministic duration taken by the
particle to reach state $1$ from state $0$ (at constant speed $\rho$).

\subsection{The resulting counting process}

Consider the maturation system with killing modelled by the $(X_t)$ process
introduced above, with maturation speed $\rho$ and killing rate
$\kappa:\dR\times[0,1]\to\dR_+$ as in \eqref{eq:def-k}. The maturation takes a
deterministic time $\tau$ given by \eqref{eq:def-tau}. A particle which begins
its maturation at time $T$ can either die with probability $1-p(T)$ where
$p(T)$ is as in \eqref{eq:def-pt-mut}, or survive with probability $p(T)$ and
achieve its maturation at time $T+\tau$. Now, consider independent particles
which begin their maturation at times $T_0,T_1,\ldots$ These particles survive
with probabilities $p(T_0),p(T_1),\ldots$, and achieve their maturation at
times $M_0=T_0+\tau,M_1=T_1+\tau,\ldots$ provided that the corresponding
$p(T_i)$ are non null. The survival durations after maturation of the mature
particles are i.i.d. random variables with common rate $\mu:\dR_+\to\dR_+$
where $\mu(t):=\kappa(t,1)$. Suppose that the initial times $T_0,T_1,\ldots$
are random and distributed according to an independent Poisson point process
on $\dR$ with intensity $\la:\dR\to\dR_+$. For any $t\geq0$, let $N_t$ be the
number of mature particles still alive at time $t$. Theorem \ref{th:maturgen},
together with \eqref{eq:def-pt-mut} and \eqref{eq:def-tau}, yields the
following result.

\begin{thm}\label{th:law}
  For any $0\leq s\leq t$ and any $m\in\dN$, $\cL(N_t\,\vert\,N_s=m) =
  \cB(m,\al(s,t))*\cP(\be(s,t))$, where
  $$
  \al(s,t):=\exp\PAR{-\int_s^t\!\mu(w)\,dw}
  $$
  and
  $$
  \be(s,t):=\int_{s-\tau}^{t-\tau}\!
  \la(u)\exp\PAR{-\int_0^{\tau}\!g(u+w,\rho w)\,dw}
  \al(u+\tau,t)\,du.
  $$
  In particular, for any $0\leq s\leq t$ and any $m\in\dN$,
  \begin{equation}\label{eq:moy}
    \moy{}{N_t\,\vert\,N_s=m}=m\al(s,t)+\be(s,t).
  \end{equation}
\end{thm}

The process $(N_t)_{t\geq0}$ is a time inhomogeneous M/M/$\infty$ queue, which
is a particular birth and death process, see
\cite{MR1996883,MR1978464,MR841197}. In an M/M/$\infty$ queue, each client is
immediately served by an independent dedicated server. When $g\equiv0$ and
$\mu$ is constant, one has $p\equiv1$, and $(N_t)_{t\geq0}$ is in that case an
M/M/$\infty$ queue with input intensity $\la$ and constant output intensity
$\mu$. When in addition $\la$ is constant, the symmetric invariant measure of
this queue is the Poisson law $\cP(\la/\mu)$. In a way, the M/M/$\infty$ queue
with constant intensities is for the Poisson process what the
Ornstein-Uhlenbeck process is for Brownian motion, see \cite[Theorem
6.14]{MR1996883}.

The role played by the M/M/$\infty$ queueing processes in our model is due to
the independence of the particles in the definition of $(N_t)_{t\geq0}$. One
can alternatively consider a non independent killing after maturation, which
can lead for example to an M/M/$1$ type queueing process. Unfortunately, the
law at fixed time of such a process, even for the time-homogeneous case, is
far more complex than the simple formula obtained in the M/M/$\infty$ case,
see \cite{MR1918702} and \cite{MR1996883}.

\begin{rem}[Negative values of $\kappa$ and input rate amplification
  heuristics]\label{re:kappa-neg}
  The expression of $\be$ above still makes sense even when $\kappa$ takes
  negative values on $[0,1)$. In that case, the quantity
  $p(u)=\exp\PAR{-\int_0^{\tau}\!g(u+w,\rho w)\,dw}$ may exceed $1$, and thus
  cannot be interpreted as a ``survival probability''. Such negative values
  can be seen in a way as a sort of ``amplification'' of the input rate
  instead of a killing, and can be incorporated into $\la$. Beware that it
  does not correspond to an immigration of particles in the counting process.
  It appears as a distortion of the input rate of the counting process.
\end{rem}

\begin{figure}[hbpt]
  \begin{center}
    \ifpdf\input{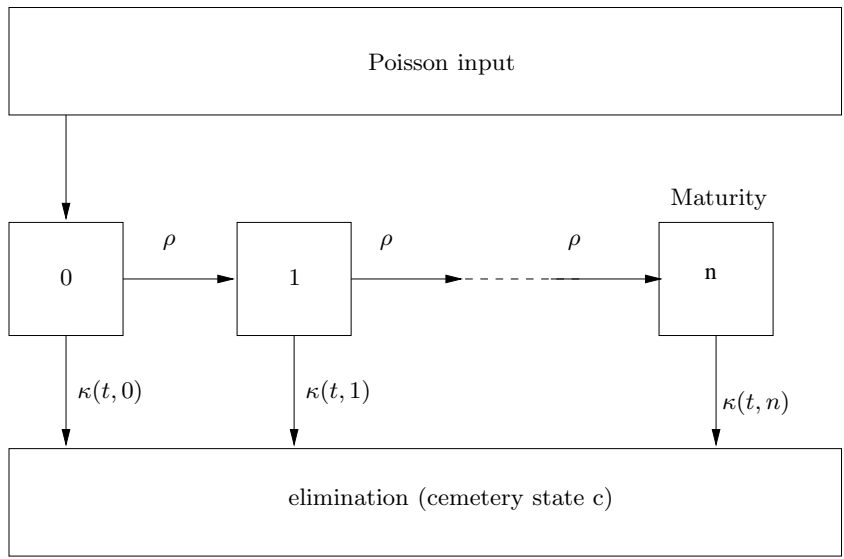}\else\input{schema1.pstex_t}\fi
    \caption{Catenary chain of compartments with Poisson input, one-way flow
      ($\rho$), and space-time-inhomogeneous killing ($\kappa$). Its limit when
      $n$ tends to infinity is given by figure \ref{fi:schem}.}
    \label{fi:mod}
  \end{center}
\end{figure}

\subsubsection{Case without killing after maturation}

Assume that $\mu\equiv0$ (no killing after maturation). In that case,
$\al(s,t)=0$ for any $0\leq s\leq t$, and the formula for $\be(s,t)$ boils
down to
\begin{equation*}
  \be(s,t) = \int_{s-\tau}^{t-\tau}\!\la(u)\exp\PAR{-\int_0^{\tau}g(u+w,\rho w)\,dw}\,du.
\end{equation*}
In that case, $t\in\dR_+\mapsto \be(s,t)$ is non decreasing since the main
integrand does not depend on $t$. This is not surprising since the particles
are never killed after maturation. Hence, on $\{N_s=0\}$, the process
$t\in\dR_+\mapsto N_t$ is non decreasing, and in particular its expected value
$\be(s,t)$ is non decreasing too. When $\la$ is constant and $g\equiv0$, we
recover a simple Poisson process of intensity $\la$.

\subsubsection{Case without killing during maturation}

Assume that $g\equiv0$ (no killing during maturation). In that case, the
formula for $\be(s,t)$ for $0\leq s\leq t$ boils down to
\begin{equation*}
  \be(s,t) = \int_{s-\tau}^{t-\tau}\!\la(u)\al(u+\tau,t)\,du.
\end{equation*}
When in addition both $\la$ and $\mu$ are constant, we recover the standard
time-homogeneous M/M/$\infty$ queue with input rate $\la$ and service rate
$\mu$, for which
\begin{equation}\label{eq:moymmi}
\al(s,t)=e^{-(t-s)\mu}
\text{\quad and\quad}
\be(s,t)=\PAR{1-e^{-(t-s)\mu}}\frac{\la}{\mu}.
\end{equation}
The Poisson measure $\cP(\la/\mu)$ is a stationary distribution of
$(N_t)_{t\geq0}$ in that case.

\subsubsection{Partial killing during maturation and constant killing after maturation}\label{ss:partial-kill}

\begin{figure}[ht] %bpt]
  \begin{center}
    \includegraphics[scale=0.6]{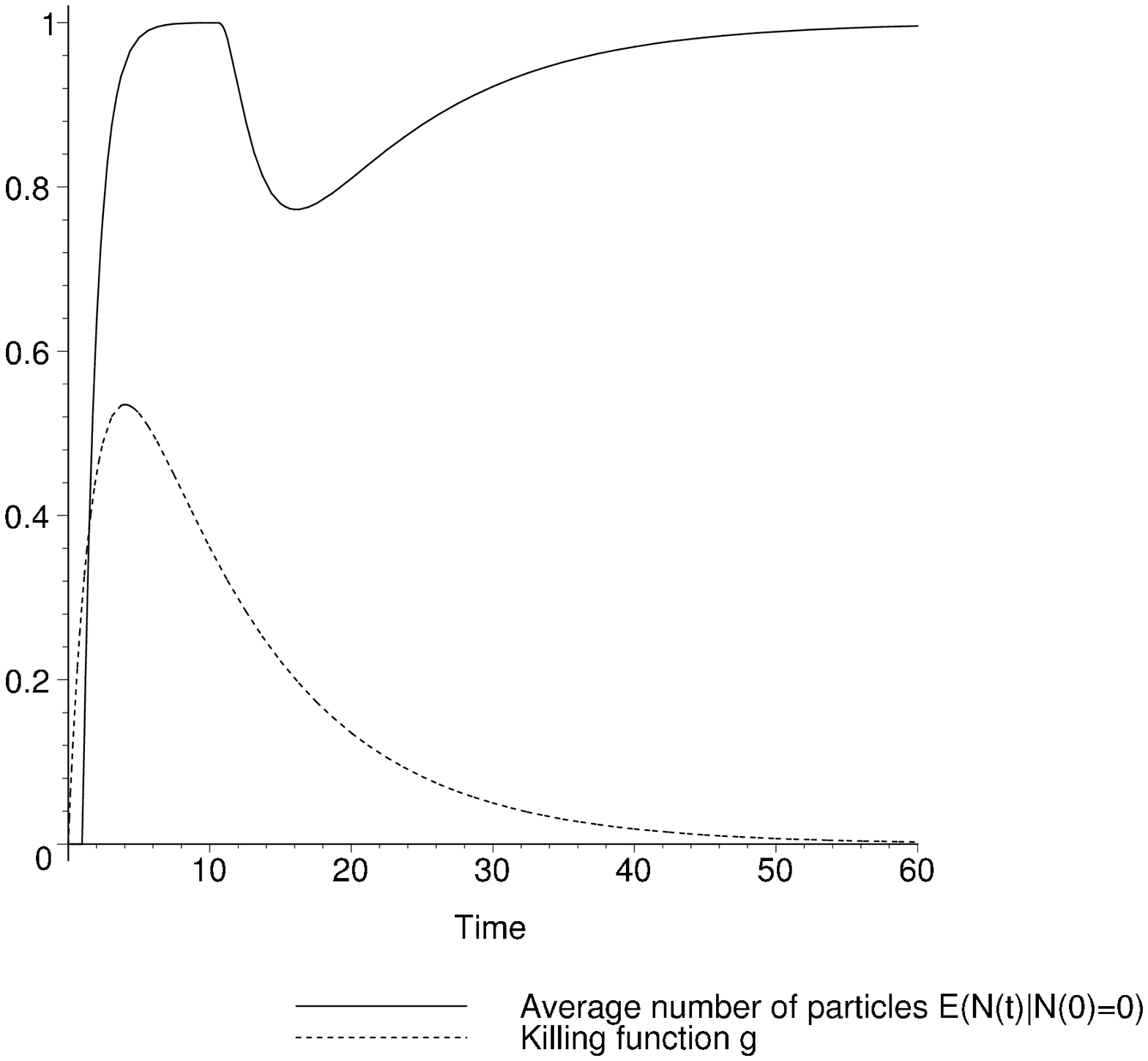}
    \caption{The plots represent the average number $\dE(N_t\,\vert\,N_0=0)$
      of Theorem \ref{th:law}, and the killing function $g$. This example
      corresponds to the case \eqref{eq:htx-gt-mu} with
      $\la(t)=\rI_{\dR_+}(t)$, $\al=\frac{1}{2}$, $\mu\equiv1$, $\rho=1$, and
      $g(t-10) =(e^{-t/10}-e^{-t/2})\rI_{\dR_+}(t)$. One can observe on the
      plot of the average number of particles three main time phases. The
      first phase corresponds to an ascendancy to a Poisson equilibrium before
      the action of the drug via $g$. The second phase corresponds to a
      decrease due to the action of the killing via function $g$ (delayed by
      the time lag is $\tau=\rho^{-1}=1$). In the third phase, the killing
      action decreases and the Poisson equilibrium is reached again.}
    \label{fi:exple}
  \end{center}
\end{figure}

Let us consider now the particular case where $\la$ is constant, and $\kappa$
is of the form
\begin{equation}\label{eq:htx-gt-mu}
  \kappa(t,x)=g(t)\rI_{[0,\de)}(x)+\mu\rI_{\{1\}}(x),
\end{equation}
where $\de\in[0,1]$, where $\mu\in\dR_+$, and where $g:\dR_+\to\dR_+$ is a
smooth function. It corresponds to a time dependent killing before the
maturation stage $\de$, and to a constant killing after maturation. No killing
is made between maturation stage $\de$ and full maturation. The formula for
$\be(s,t)$ when $0\leq s\leq t$ boils down to
\begin{equation}\label{eq:rt-gt-mu}
  \be(s,t) = \la \int_{s-\tau}^{t-\tau}\!%
  \exp\PAR{-\mu(t-\tau-u)-\int_0^{\de\tau}\!\!\!g(u+w)\,dw}\,du.
\end{equation}
When $g\equiv0$, this formula reduces to the classical M/M/$\infty$ average
queue length \eqref{eq:moymmi}. Assume instead that function $g$ in
\eqref{eq:htx-gt-mu} only vanishes at infinity. Then the two formulas
\eqref{eq:rt-gt-mu} and \eqref{eq:moymmi} for $\be$ are equivalent when $t$
goes to $+\infty$. In particular, the Poisson measure $\cP(\la/\mu)$ is a
stationary distribution of $(N_t)_{t\geq0}$ in that case. One can see on
figure \ref{fi:exple} that this Poisson equilibrium is quickly reached. It is
possible to quantify the speed of convergence in total variation norm or in
entropy, see \cite{MR2247924,MR1996883}.

\subsubsection{Identifiability and invariance by some transformations}
\label{ss:identif}

The dynamics of $N:=(N_t)_{t\geq0}$ is fully described by the quadruple
$(\tau,\mu,g,\la)$, where $\tau:=1/\rho$. It is quite natural to ask about the
injectivity of the map
$$
(\tau,\mu,g,\la)\mapsto\cL(N\,\vert\,N_0).
$$
According to Theorem \ref{th:law}, the law $\cL(N\,\vert\,N_0)$ is completely
prescribed by the triple $(\tau,\mu,\la_g)$, where $\la_g:\dR\to\dR_+$ is
defined by
$$
\la_g(t):=\la(t)p(t)=\la(t)\exp\PAR{-\int_0^\tau\!g(t+w,\rho w)\,dw}.
$$
Thus, the couple $(\la,g)$ cannot be identified, since it acts on the dynamics
via the compound $\la_g$. The action of $g$ can be compensated by $\la$ and
vice versa. Namely, suppose that $g$ can be written $g=g_1+g_2$, where $g_1$
and $g_2$ are non-negative functions. Then, the two models corresponding
respectively to $(\tau,\mu,g,\la)$ and $(\tau,\mu,g_2,\la_{g_1})$ are
indistinguishable. The extreme case corresponds to $(\tau,\mu,0,\la_g)$, for
which the entire killing function $g$ is merged into the input rate function
$\la$. Let us analyze some special cases. For any $\te>1$, let us consider the
transformation which replaces $(\la,g)$ by $(\la^\te,g^\te)$ defined by
\begin{equation*}
  \la^\te:=\te\la
  \text{\quad and \quad}
  g^\te:=g+\rho\,\log(\te).
\end{equation*}
Function $\mu$ and parameter $\rho$ are left unchanged, and one can check that
$\la^\te_{g^\te}=\la_g$. Therefore, the dynamics is invariant by this
transformation: the models corresponding to $(\tau,\mu,g,\la)$ and
$(\tau,\mu,g^\te,\la^\te)$ are indistinguishable, despite their distinct
physiological meanings. A multiplicative perturbation of the input intensity
$\la$ corresponds to an additive perturbation of the killing function $g$
during maturation. Hence, one can decide to normalize the parametrization by
taking for example $\la\equiv1$. Namely, $(\tau,\mu,g,\lambda)$ is
indistinguishable from $(\tau,\mu,G,1)$ where
$$
G(v,y)%
=g(v,y)%
-\frac{1}{\tau}\log\lambda\left(v-\frac{y}{\rho}\right).
$$

\section{A pharmacokinetics/pharmacodynamics example in cancerology}
\label{se:example}

Catenary chains of compartments, as depicted in figure \ref{fi:family2}, are
widely used in the literature by kineticists for the modelling of anticancer
drug toxicity, see for instance
\cite{dconc,dayneka-jusko,friberg-karlsson-al,karlsson,Friberg2,gobburu-jusko}
and references therein. Most anticancer drugs have toxic effects on white
blood cells (myelosuppression). Neutrophils are particular white blood cells
used as toxicity markers. The catenary chain of compartments models the
maturation process of neutrophils in the bone marrow. Each maturation stage
corresponds to a specific position in the chain. The killing rate during
maturation corresponds to the drug toxicity on neutrophils. The last
compartment of the chain is the only observed (blood), and the others
correspond to hidden positions (bone marrow). In practice, the input and
transit rates in the catenary chain of compartments are unknown parameters.
The goal of the kineticist is thus to control the content of the last
compartment by his action on the anticancer drug dosage regimen. Consider for
simplicity that both the production rate $\la$, the transit rate $\rho$, and
the elimination rate in blood $\mu$ are constant. The killing rate $\kappa$
depends typically on the time profile $q:\dR_+\to\dR$ of the anticancer drug
concentration in blood (i.e. kinetics). The first step of the
Pharmacokinetics/Pharmacodynamics study is to find, for a given kinetics, the
best values of the parameters in view of the data (in blood). This first step
can be addressed in several ways.

Approach (I). This traditional approach involves a deterministic finite
catenary chain of $n$ compartments, as presented in example
\ref{xp:finite-catenary} page \pageref{xp:finite-catenary}.
% figure \ref{fi:family2} gives the stochastic interpretation
It is customary to model the toxic effect of the anticancer drug on
neutrophils via a killing rate $\kappa(t,i) = \ga q(t)$ for $i\leq n_{0}<n$,
$\kappa(t,n) = \mu $, and $0$ elsewhere. The vector
$Q(t):=(Q_1(t),\ldots,Q_n(t))$, which represents the number of neutrophils in
each compartments of the chain, solves a system of $n$ ordinary differential
equations (ODE). The available data concerns only the last compartment
(blood), labelled $n$. The function $t\mapsto Q_n(t)$ depends on the modelling
parameters $\la$, $\rho$, $\mu$, $\ga$, $n_0$. The main drawback of this
approach is that the number of compartments $n$ is unknown, and that the
system of equation is heavy.

Approach (II). In this approach, we circumvent the problem of the choice of
$n$ by considering the continuous limit in $n$, as presented in Section
\ref{ss:stoch-interp}. This limit procedure leads to killing rates $\kappa$ as
in equation \eqref{eq:def-k}. More precisely, we mimic the cutoff used in the
approach (I) by considering equation \eqref{eq:htx-gt-mu} with $g(t,x)=\ga
q(t)\rI_{[0,\de]}(x)$ where $\ga>0$ and where $\de\in[0,1)$ are constants. The
continuous limit in $n$ is constructed at the level of the stochastic
interpretation. As a consequence, the relevant quantity is the mean of the
counting process of Theorem \ref{th:law}, given by equation \eqref{eq:moy}.
The initial value $m$ in \eqref{eq:moy} must be averaged with respect to the
equilibrium without drug, i.e. the Poisson measure $\cP(\la/\mu)$. The
quantity $Q_n(t)$ of the approach (I) is thus replaced now by the following
\emph{explicit} alternative quantity:
\begin{equation}\label{eq:Q-antic}
  Q_\infty(t):=\al(0,t)+\be(0,t)\frac{\la}{\mu},
\end{equation}
where the functions $\al$ and $\be$ are as in Theorem \ref{th:law}.

\begin{rem}
  Actually, one can use a Hill-type transform (see for instance
  \cite{MR1908418}) to model the dependency of the killing rate with respect
  to the drug kinetics. However, for low drug concentration values, it is
  customary to simplify the approach by considering a linear dependency with a
  possible cutoff in space, as in the approaches (I) and (II) above. We
  emphasize the fact that even if the instantaneous killing rate depends
  linearly on the drug kinetics, the global killing effect is nonlinear with
  respect to the kinetics.
\end{rem}

We compare below, on a simple example, the two approaches defined above. The
approach (I) leads to the resolution of a system of ODE, whereas the approach
(II) leads to the computation of explicit integral transforms given by
\eqref{eq:Q-antic}.

\begin{figure}[hbpt]
  \begin{center}
    \includegraphics[scale=0.5,angle=-90]{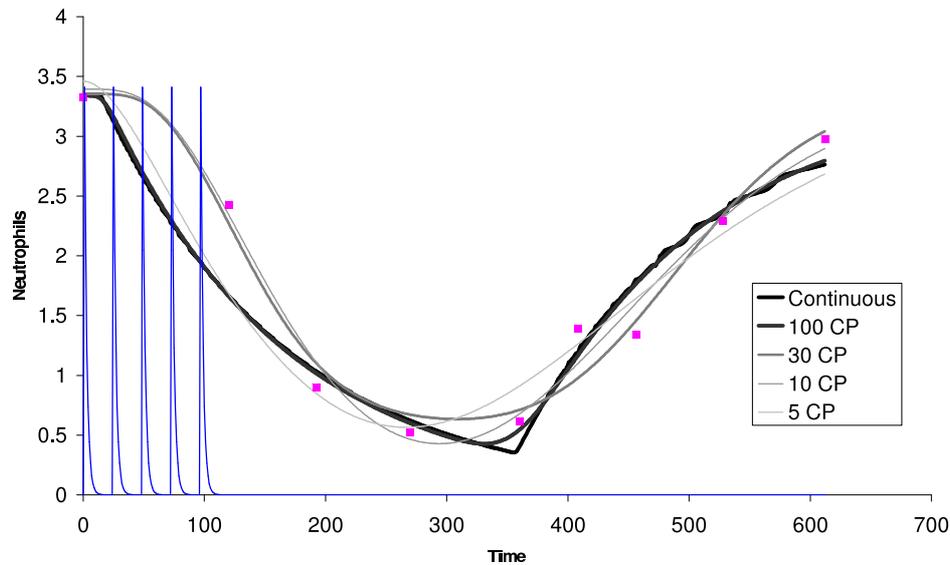}
  \end{center}
  \caption{The curve with five peaks is the drug concentration curve $t\mapsto
    q(t)$ (kinetics). The neutrophils counts observed in a patient are
    represented by squares. The mean neutrophils counts predicted with models
    with 5, 10, 30, 100, $\infty$ compartments are represented with different
    lines. On the whole, models with at least ten compartments give a
    reasonable description of observed points. We emphasize that in general,
    the adequate number of compartments is not known a priori, and can vary
    from one patient to other. The estimation of an adequate number of
    compartments from the data has a cost. The continuous model avoids this
    problem.}
  \label{fi:courbescomp}
\end{figure}

In order to mimic practical situations, we used a synthetic dataset obtained
by perturbing a dataset extracted from \cite{dconc} (we did not have the
permission to use directly the dataset of \cite{dconc}). In this study, an
anticancer drug was administrated by 30-min intravenous infusions to patients,
on five consecutive days, see figure \ref{fi:courbescomp}. For a typical
patient, the drug concentration at time $t$ was proportional to :
$$
q(t)=\sum_{d=0}^{4} \PAR{1-e^{-\alpha(t-24d)}}\rI_{24d+[0,1/2]}(t) %
+\PAR{1-e^{-\alpha/2}}e^{-\beta(t-24d-1/2)}\rI_{24d+[1/2,\infty]}(t),
$$
with $\alpha= 1.86$ and $\beta=0.51$. For simplicity, we restrict our analysis
in the sequel to a single patient. We first analyzed this dataset using
approach (I) and approach (II), as described above. Let us denote by $Y_j$ the
observed neutrophils counts in blood at time $t_j$.

Approach (I). The vector valued function $t\mapsto(Q_1(t),\ldots,Q_n(t))$
solves the following system of ordinary differential equations:
$$
\left\{
\begin{array}{lll}
\pd_t Q_{1}(t) &
=& \la -\rho Q_{1}(t)-\ga q(t) Q_{1}(t)\\
\pd_t Q_{2}(t) &
=& \rho Q_{1}(t)-\rho Q_{2}(t)-\ga q(t) Q_{2}(t)\\
&\vdots&\\
\pd_t Q_{n_{0}}(t) &
=& \rho Q_{n_{0}-1}(t)-\rho Q_{n_{0}}(t)-\ga q(t) Q_{n_{0}}(t)\\
&\vdots&\\
\pd_t Q_{n}(t) &
=& \rho Q_{n-1}(t)-\mu Q_{n}(t)
\end{array}
\right.
$$
At time $t=0$, the drug concentration $q$ is null, and the system is at
equilibrium. Consequently, $\pd_{t=0} Q_{i}=0$ for $i\in\{1,\ldots,n\}$. The
initial condition is thus $Q_{1}(0)=\cdots=Q_{n-1}(0)=\la/\rho$, and
$Q_{n}(0)=\la/\mu$. These systems were numerically integrated using the
\textsc{Fortran} package \textsc{ODEpack}, see for instance
\cite{Radhakrishnan,MR689694,MR751604}. Models with n= 5, 10, 30 and 100
compartments were used. For each $n$-compartments model, the parameters
($\la$, $\rho$, $\mu$, $\ga$, $n_{0}$) were estimated with ordinary
least-squares :
$$
\arg\inf_{\la,\rho,\mu,\ga,n_{0}}\sum_{j}\PAR{Y_{j}-Q_{n}(t_{j})}^{2}.
$$
Since $n_{0}$ is an integer, all possible values of $n_{0}$, from $1$ to
$n-1$, were screened and the value of $n_{0}$ giving the minimum residual sum
of squares was then selected.

Approach (II). The explicit expression of \eqref{eq:Q-antic} provides the
following formula for the average number of neutrophils in blood at time $t$:
$$
Q_\infty(t) = \frac{\la}{\mu} e^{-\mu t} +\la e^{\mu(t-\tau)}
\int_{-\tau}^{t-\tau}\!\exp\PAR{\mu u - \ga\int_{0}^{\de\tau}q(u+w)dw}\,du,
$$
where $\tau=1/\rho$. Since $q$ has an easily computable primitive, the inner
integral in the expression above is explicit. The outer integral was
numerically evaluated using the Clenshaw-Curtis quadrature method. The
parameters ($\la$, $\rho$, $\mu$, $\ga$, $\de$) were estimated with ordinary
least-squares. Note that the continuous parameter $\de\in[0,1)$ plays here the
role of the discrete parameter $n_0< n$ of the approach (I).

The curves obtained with the different models are represented in figure
\ref{fi:courbescomp}. On the whole, the different models approximately give
the same curves. Despite its appearance, $t\mapsto Q_\infty(t)$ is smooth. The
curves provided by the approach (I) with 100 compartments in one hand and the
approach (II) in the other hand are nearly the same. The obtained residuals
sum of squares are given in table \ref{table:tempscalcul}. The curves are
given by figure \ref{fi:courbescomp}. The best fit was obtained with the
10-compartments model. Models with more compartments, including the $\infty$
compartments model of approach (II) give approximately the same quality of
fit. The number of parameter to estimate is the same for both approach (I) and
(II). However, the estimation of $n_{0}$ in approach (I) leads to a large
amount of extra ODE integration. The amounts of time needed for the ODE
integration are at least twice as long for models with more than 10
compartments than for the integral evaluation of approach (II). Keeping in
mind that such studies usually include several tens of patients, the approach
(II) that we proposed offers a reasonable alternative to compartmental models.

\begin{table}[h]
\begin{center}
\begin{tabular}{llllll}
\hline
{}&5 CP&10 CP&30 CP&100 CP&Continuous\\
\hline
Computation time& 0.39  &  0.99   & 11.803 & 235.12 &0.51\\
Sum of Squares &  0.99 & 0.25 & 0.32& 0.32 & 0.33 \\
\hline
\end{tabular}
\end{center}
\caption{The duration of the ODE/integral evaluation. For compartmental models,
  these durations have been obtained as $(n-1)\times$ duration for a single
  evaluation of the ODE to take into account the screening needed to estimate
  $n_{0}$. The obtained sum-of-squares are given in the second row.} 
\label{table:tempscalcul}
\end{table}

\section*{Acknowledgements}

The authors sincerely thank Prof. Laurent Miclo for stimulating and helpful
discussions on inhomogeneous jump processes and related topics at the early
development of this work. The article greatly benefited from the comments of
two anonymous reviewers.

%
%\nocite{*}
\addcontentsline{toc}{section}{\refname}%
{
\footnotesize
%\bibliography{\jobname}
\bibliography{neutro}
}
\bibliographystyle{amsplain}

\begin{center}
  \hrule
\end{center}

{\small
  \noindent Djalil \textsc{Chafa\"\i} and Didier \textsc{Concordet}\\
  \textbf{Address:} UMR 181 INRA/ENVT Physiopathologie et Toxicologie
  Exp\'erimentales, \\
  \'Ecole Nationale V\'et\'erinaire de Toulouse, \\
  23 Chemin des Capelles, F-31076, Toulouse \textsc{Cedex} 3, France.\\
  \textbf{E-mail:} \url{mailto:d.chafai@envt.fr.nospam}\\
  \textbf{Address:} UMR 5583 CNRS/UPS Laboratoire de Statistique et
  Probabilit\'es, \\
  Institut de Math\'ematiques de Toulouse, Universit\'e Paul Sabatier, \\
  118 route de Narbonne, F-31062, Toulouse, \textsc{Cedex} 4, France.\\
  \textbf{E-mail:} \url{mailto:chafai(AT)math.univ-toulouse.fr}\\
  \textbf{Web:} \url{http://www.math.univ-toulouse/~chafai}\\
}%small

\end{document}

%% file: macros.tex
\newcommand{\EM}{\ensuremath}

\newcommand{\dE}{\EM{\mathbb{E}}}

\newcommand{\dN}{\EM{\mathbb{N}}}

\newcommand{\dP}{\EM{\mathbb{P}}}

\newcommand{\dR}{\EM{\mathbb{R}}}

\newcommand{\rI}{\EM{\mathrm{I}}}

\newcommand{\cB}{\EM{\mathcal{B}}}

\newcommand{\cL}{\EM{\mathcal{L}}}

\newcommand{\cP}{\EM{\mathcal{P}}}

\newcommand{\cS}{\EM{\mathcal{S}}}

\newcommand{\bE}{\EM{\mathbf{E}}}

\newcommand{\bI}{\EM{\mathbf{I}}}

\newcommand{\bL}{\EM{\mathbf{L}}}
\newcommand{\bM}{\EM{\mathbf{M}}}

\newcommand{\bR}{\EM{\mathbf{R}}}

\newcommand{\al}{\alpha}
\newcommand{\be}{\beta}

\newcommand{\de}{\delta}
\newcommand{\ga}{\gamma}

\newcommand{\la}{\lambda}

\newcommand{\te}{\theta}

\newcommand{\veps}{\varepsilon}

\newcommand{\p}[4]{{#3}\!\left#1{#4}\right#2} 

\newcommand{\PENT}[1]{\EM{\lfloor#1\rfloor}} % partie entiere
\newcommand{\ABS}[1]{\EM{{\left| #1 \right|}}} % |1|
\newcommand{\BRA}[1]{\EM{{\left\{#1\right\}}}} % {1}
\newcommand{\PAR}[1]{\EM{{\left(#1\right)}}} % (1)
\newcommand{\pd}{\EM{\partial}} % derivee partielle
\newcommand{\SBRA}[1]{\EM{{\left[#1\right]}}} % [1]
\newcommand{\moyf}[1]{\bE_{#1}}
\newcommand{\moy}[2]{\p(){\moyf{#1}}{#2}}

\newcommand{\Det}[1]{\mathrm{Det}\,}

\renewcommand{\leq}{\leqslant}
\renewcommand{\geq}{\geqslant}

%% file: compgen.pstex_t
\begin{picture}(0,0)%
\includegraphics{compgen.pstex}%
\end{picture}%
\setlength{\unitlength}{4144sp}%
\begingroup\makeatletter\ifx\SetFigFont\undefined%
\gdef\SetFigFont#1#2#3#4#5{%
  \reset@font\fontsize{#1}{#2pt}%
  \fontfamily{#3}\fontseries{#4}\fontshape{#5}%
  \selectfont}%
\fi\endgroup%
\begin{picture}(5562,2679)(2146,-4348)
\put(6391,-3886){\makebox(0,0)[lb]{\smash{{\SetFigFont{12}{14.4}{\familydefault}{\mddefault}{\updefault}{\color[rgb]{0,0,0}$\kappa_i(t)$}%
}}}}
\put(6931,-2761){\makebox(0,0)[lb]{\smash{{\SetFigFont{12}{14.4}{\familydefault}{\mddefault}{\updefault}{\color[rgb]{0,0,0}$\rho_{i,j}(t)$}%
}}}}
\put(6391,-2086){\makebox(0,0)[lb]{\smash{{\SetFigFont{12}{14.4}{\familydefault}{\mddefault}{\updefault}{\color[rgb]{0,0,0}$\la_i(t)$}%
}}}}
\put(6931,-3391){\makebox(0,0)[lb]{\smash{{\SetFigFont{12}{14.4}{\familydefault}{\mddefault}{\updefault}{\color[rgb]{0,0,0}$\rho_{j,i}(t)$}%
}}}}
\put(3376,-3031){\makebox(0,0)[lb]{\smash{{\SetFigFont{12}{14.4}{\familydefault}{\mddefault}{\updefault}{\color[rgb]{0,0,0}$Q_i(t)$}%
}}}}
\put(6121,-3076){\makebox(0,0)[lb]{\smash{{\SetFigFont{12}{14.4}{\familydefault}{\mddefault}{\updefault}{\color[rgb]{0,0,0}$Q_i(t)$}%
}}}}
\put(2161,-2761){\makebox(0,0)[lb]{\smash{{\SetFigFont{12}{14.4}{\familydefault}{\mddefault}{\updefault}{\color[rgb]{0,0,0}$Q_j(t)\rho_{j,i}(t)$}%
}}}}
\put(2161,-3391){\makebox(0,0)[lb]{\smash{{\SetFigFont{12}{14.4}{\familydefault}{\mddefault}{\updefault}{\color[rgb]{0,0,0}$Q_i(t)\rho_{i,j}(t)$}%
}}}}
\put(3691,-3931){\makebox(0,0)[lb]{\smash{{\SetFigFont{12}{14.4}{\familydefault}{\mddefault}{\updefault}{\color[rgb]{0,0,0}$Q_i(t)\kappa_i(t)$}%
}}}}
\put(3646,-2086){\makebox(0,0)[lb]{\smash{{\SetFigFont{12}{14.4}{\familydefault}{\mddefault}{\updefault}{\color[rgb]{0,0,0}$\la_i(t)$}%
}}}}
\end{picture}%

%% file: fincat.pstex_t
\begin{picture}(0,0)%
\includegraphics{fincat.pstex}%
\end{picture}%
\setlength{\unitlength}{4144sp}%
\begingroup\makeatletter\ifx\SetFigFont\undefined%
\gdef\SetFigFont#1#2#3#4#5{%
  \reset@font\fontsize{#1}{#2pt}%
  \fontfamily{#3}\fontseries{#4}\fontshape{#5}%
  \selectfont}%
\fi\endgroup%
\begin{picture}(6774,1104)(2149,-3493)
\put(8011,-3256){\makebox(0,0)[lb]{\smash{{\SetFigFont{12}{14.4}{\rmdefault}{\mddefault}{\updefault}{\color[rgb]{0,0,0}$\kappa_n(t)$}%
}}}}
\put(3511,-2581){\makebox(0,0)[lb]{\smash{{\SetFigFont{12}{14.4}{\rmdefault}{\mddefault}{\updefault}{\color[rgb]{0,0,0}$\rho_1(t)$}%
}}}}
\put(6391,-2581){\makebox(0,0)[lb]{\smash{{\SetFigFont{12}{14.4}{\rmdefault}{\mddefault}{\updefault}{\color[rgb]{0,0,0}$\rho_{i}(t)$}%
}}}}
\put(4816,-2581){\makebox(0,0)[lb]{\smash{{\SetFigFont{12}{14.4}{\rmdefault}{\mddefault}{\updefault}{\color[rgb]{0,0,0}$\rho_{i-1}(t)$}%
}}}}
\put(5266,-3301){\makebox(0,0)[lb]{\smash{{\SetFigFont{12}{14.4}{\rmdefault}{\mddefault}{\updefault}{\color[rgb]{0,0,0}$\kappa_i(t)$}%
}}}}
\put(2341,-3301){\makebox(0,0)[lb]{\smash{{\SetFigFont{12}{14.4}{\rmdefault}{\mddefault}{\updefault}{\color[rgb]{0,0,0}$\kappa_1(t)$}%
}}}}
\put(7876,-2626){\makebox(0,0)[lb]{\smash{{\SetFigFont{12}{14.4}{\rmdefault}{\mddefault}{\updefault}{\color[rgb]{0,0,0}$\rho_n(t)$}%
}}}}
\put(2206,-2581){\makebox(0,0)[lb]{\smash{{\SetFigFont{12}{14.4}{\rmdefault}{\mddefault}{\updefault}{\color[rgb]{0,0,0}$\lambda$}%
}}}}
\put(2836,-2806){\makebox(0,0)[lb]{\smash{{\SetFigFont{12}{14.4}{\rmdefault}{\mddefault}{\updefault}{\color[rgb]{0,0,0}$Q_1(t)$}%
}}}}
\put(8416,-2806){\makebox(0,0)[lb]{\smash{{\SetFigFont{12}{14.4}{\rmdefault}{\mddefault}{\updefault}{\color[rgb]{0,0,0}$Q_n(t)$}%
}}}}
\put(5671,-2806){\makebox(0,0)[lb]{\smash{{\SetFigFont{12}{14.4}{\rmdefault}{\mddefault}{\updefault}{\color[rgb]{0,0,0}$Q_i(t)$}%
}}}}
\end{picture}%

%% file: stoch.pstex_t
\begin{picture}(0,0)%
\includegraphics{stoch.pstex}%
\end{picture}%
\setlength{\unitlength}{4144sp}%
\begingroup\makeatletter\ifx\SetFigFont\undefined%
\gdef\SetFigFont#1#2#3#4#5{%
  \reset@font\fontsize{#1}{#2pt}%
  \fontfamily{#3}\fontseries{#4}\fontshape{#5}%
  \selectfont}%
\fi\endgroup%
\begin{picture}(5064,1950)(1159,-4090)
\put(1441,-2311){\makebox(0,0)[lb]{\smash{{\SetFigFont{12}{14.4}{\rmdefault}{\mddefault}{\updefault}{\color[rgb]{0,0,0}M/M/$\infty$ queue}%
}}}}
\put(4501,-2311){\makebox(0,0)[lb]{\smash{{\SetFigFont{12}{14.4}{\rmdefault}{\mddefault}{\updefault}M/M/$\infty$ queue}}}}
\put(3061,-3616){\makebox(0,0)[lb]{\smash{{\SetFigFont{12}{14.4}{\rmdefault}{\mddefault}{\updefault}{\color[rgb]{0,0,0}Random walks}%
}}}}
\put(2791,-4021){\makebox(0,0)[lb]{\smash{{\SetFigFont{12}{14.4}{\rmdefault}{\mddefault}{\updefault}{\color[rgb]{0,0,0}Graph of compartments}%
}}}}
\end{picture}%

%% file: schema2.pstex_t
\begin{picture}(0,0)%
\includegraphics{schema2.pstex}%
\end{picture}%
\setlength{\unitlength}{4144sp}%
\begingroup\makeatletter\ifx\SetFigFont\undefined%
\gdef\SetFigFont#1#2#3#4#5{%
  \reset@font\fontsize{#1}{#2pt}%
  \fontfamily{#3}\fontseries{#4}\fontshape{#5}%
  \selectfont}%
\fi\endgroup%
\begin{picture}(6504,1705)(529,-2969)
\put(1261,-1591){\makebox(0,0)[lb]{\smash{{\SetFigFont{12}{14.4}{\rmdefault}{\mddefault}{\updefault}{\color[rgb]{0,0,0}Input}%
}}}}
\put(4276,-1636){\makebox(0,0)[lb]{\smash{{\SetFigFont{12}{14.4}{\rmdefault}{\mddefault}{\updefault}{\color[rgb]{0,0,0}Output}%
}}}}
\put(5761,-1816){\makebox(0,0)[lb]{\smash{{\SetFigFont{12}{14.4}{\rmdefault}{\mddefault}{\updefault}{\color[rgb]{0,0,0}Maturity}%
}}}}
\put(1171,-2086){\makebox(0,0)[lb]{\smash{{\SetFigFont{12}{14.4}{\rmdefault}{\mddefault}{\updefault}{\color[rgb]{0,0,0}time $T$}%
}}}}
\put(4051,-2041){\makebox(0,0)[lb]{\smash{{\SetFigFont{12}{14.4}{\rmdefault}{\mddefault}{\updefault}{\color[rgb]{0,0,0}time $T+\tau$}%
}}}}
\put(1711,-2896){\makebox(0,0)[lb]{\smash{{\SetFigFont{12}{14.4}{\rmdefault}{\mddefault}{\updefault}{\color[rgb]{0,0,0}Elimination with probability $1-p(T)$}%
}}}}
\put(5131,-2896){\makebox(0,0)[lb]{\smash{{\SetFigFont{12}{14.4}{\rmdefault}{\mddefault}{\updefault}{\color[rgb]{0,0,0}Elimination with rate $\mu$}%
}}}}
\put(2746,-1816){\makebox(0,0)[lb]{\smash{{\SetFigFont{12}{14.4}{\rmdefault}{\mddefault}{\updefault}{\color[rgb]{0,0,0}Lag $\tau$}%
}}}}
\end{picture}%

%% file: schema0.pstex_t
\begin{picture}(0,0)%
\includegraphics{schema0.pstex}%
\end{picture}%
\setlength{\unitlength}{4144sp}%
\begingroup\makeatletter\ifx\SetFigFont\undefined%
\gdef\SetFigFont#1#2#3#4#5{%
  \reset@font\fontsize{#1}{#2pt}%
  \fontfamily{#3}\fontseries{#4}\fontshape{#5}%
  \selectfont}%
\fi\endgroup%
\begin{picture}(6812,1431)(2146,-3493)
\put(7201,-3256){\makebox(0,0)[lb]{\smash{{\SetFigFont{12}{14.4}{\rmdefault}{\mddefault}{\updefault}{\color[rgb]{0,0,0}$\mu(t)=\kappa^{(n)}(t,n)$}%
}}}}
\put(5041,-3301){\makebox(0,0)[lb]{\smash{{\SetFigFont{12}{14.4}{\rmdefault}{\mddefault}{\updefault}{\color[rgb]{0,0,0}$\kappa^{(n)}(t,i)$}%
}}}}
\put(7876,-2626){\makebox(0,0)[lb]{\smash{{\SetFigFont{12}{14.4}{\rmdefault}{\mddefault}{\updefault}{\color[rgb]{0,0,0}$\rho^{(n)}$}%
}}}}
\put(2161,-3301){\makebox(0,0)[lb]{\smash{{\SetFigFont{12}{14.4}{\rmdefault}{\mddefault}{\updefault}{\color[rgb]{0,0,0}$\kappa^{(n)}(t,0)$}%
}}}}
\put(3511,-2581){\makebox(0,0)[lb]{\smash{{\SetFigFont{12}{14.4}{\rmdefault}{\mddefault}{\updefault}{\color[rgb]{0,0,0}$\rho^{(n)}$}%
}}}}
\put(5131,-2581){\makebox(0,0)[lb]{\smash{{\SetFigFont{12}{14.4}{\rmdefault}{\mddefault}{\updefault}{\color[rgb]{0,0,0}$\rho^{(n)}$}%
}}}}
\put(6391,-2581){\makebox(0,0)[lb]{\smash{{\SetFigFont{12}{14.4}{\rmdefault}{\mddefault}{\updefault}{\color[rgb]{0,0,0}$\rho^{(n)}$}%
}}}}
\put(8236,-2221){\makebox(0,0)[lb]{\smash{{\SetFigFont{12}{14.4}{\rmdefault}{\mddefault}{\updefault}{\color[rgb]{0,0,0}maturity}%
}}}}
\end{picture}%

%% file: schema1.pstex_t
\begin{picture}(0,0)%
\includegraphics{schema1.pstex}%
\end{picture}%
\setlength{\unitlength}{3108sp}%
\begingroup\makeatletter\ifx\SetFigFont\undefined%
\gdef\SetFigFont#1#2#3#4#5{%
  \reset@font\fontsize{#1}{#2pt}%
  \fontfamily{#3}\fontseries{#4}\fontshape{#5}%
  \selectfont}%
\fi\endgroup%
\begin{picture}(6594,4389)(1339,-5203)
\put(2566,-2716){\makebox(0,0)[lb]{\smash{{\SetFigFont{9}{10.8}{\familydefault}{\mddefault}{\updefault}{\color[rgb]{0,0,0}$\rho$}%
}}}}
\put(4276,-2716){\makebox(0,0)[lb]{\smash{{\SetFigFont{9}{10.8}{\familydefault}{\mddefault}{\updefault}{\color[rgb]{0,0,0}$\rho$}%
}}}}
\put(6976,-4021){\makebox(0,0)[lb]{\smash{{\SetFigFont{9}{10.8}{\familydefault}{\mddefault}{\updefault}{\color[rgb]{0,0,0}$\kappa(t,n)$}%
}}}}
\put(3556,-3031){\makebox(0,0)[lb]{\smash{{\SetFigFont{9}{10.8}{\familydefault}{\mddefault}{\updefault}{\color[rgb]{0,0,0}1}%
}}}}
\put(1756,-3031){\makebox(0,0)[lb]{\smash{{\SetFigFont{9}{10.8}{\familydefault}{\mddefault}{\updefault}{\color[rgb]{0,0,0}0}%
}}}}
\put(3691,-3931){\makebox(0,0)[lb]{\smash{{\SetFigFont{9}{10.8}{\familydefault}{\mddefault}{\updefault}{\color[rgb]{0,0,0}$\kappa(t,1)$}%
}}}}
\put(1891,-3931){\makebox(0,0)[lb]{\smash{{\SetFigFont{9}{10.8}{\familydefault}{\mddefault}{\updefault}{\color[rgb]{0,0,0}$\kappa(t,0)$}%
}}}}
\put(5761,-2716){\makebox(0,0)[lb]{\smash{{\SetFigFont{9}{10.8}{\familydefault}{\mddefault}{\updefault}{\color[rgb]{0,0,0}$\rho$}%
}}}}
\put(3556,-4786){\makebox(0,0)[lb]{\smash{{\SetFigFont{9}{10.8}{\rmdefault}{\mddefault}{\updefault}{\color[rgb]{0,0,0}elimination (cemetery state c)}%
}}}}
\put(4186,-1321){\makebox(0,0)[lb]{\smash{{\SetFigFont{9}{10.8}{\rmdefault}{\mddefault}{\updefault}{\color[rgb]{0,0,0}Poisson input}%
}}}}
\put(6571,-2401){\makebox(0,0)[lb]{\smash{{\SetFigFont{9}{10.8}{\rmdefault}{\mddefault}{\updefault}{\color[rgb]{0,0,0}Maturity}%
}}}}
\end{picture}%